\documentclass{amsart}

\usepackage{srcltx}
\usepackage{amsfonts,euscript,eufrak,latexsym,amssymb,amsbsy,amsthm,stmaryrd}
 \usepackage{rotating}
\usepackage{graphicx}
 \usepackage{url}
  \usepackage{tikz}
  \usepackage{hyperref}

\newcommand{\B}{\mathcal{B}}

\begin{document}

\title{
Finite Element Methods for Elastic Contact: Penalty and Nitsche
}
\author{Tom Gustafsson,  Rolf Stenberg}

\begin{abstract}
We consider two methods for treating elastic contact problems with the finite element method; the penalty method and Nitsche's method. For the penalty method we discuss how the penalty parameter should be chosen. Both the theoretical analysis and numerical examples show that an optimal convergence rate cannot be achieved. The method is contrasted to that of Nitsche which is optimally convergent. We also give the derivation of Nitsche's method by a very simple consistency correction of the penalty method.
  \end{abstract}


\maketitle

\section{Introduction}

 The paper {\em Variational methods for the solution of problems of equilibrium and vibrations} \cite{Courant}, presented by Richard Courant at the American Mathematics Society meeting in Washington D.C. in 1941, is widely considered as a starting point in two fields of applied mathematics, finite element methods and the use of the penalty method in optimization. 

For the case of optimization this appears to be the case, see, e.g.,  the classical books \cite{Gill, Polak}.

In the field of finite elements, the continuous piecewise linear basis functions are often referred to the "Courant triangle" or "Courant hat function". The paper by Courant came out in print in 1943,  and to celebrate the 50th anniversary, a conference was organized in Jyv\"askyl\"a \cite{KNS}. Ivo Babu\v{s}ka gave an invited talk, in which he outlined the history of the finite element method. Surprisingly, he discovered that none of the earliest paper on the mathematical foundations of the FE cite Courants paper. Furthermore, he pointed out that piecewise linear approximation method was not presented at the conference, it was introduced as an appendix in the final paper. 

The same comment is given by Gilbert Strang in the  first textbook (Strang \& Fix \cite{StrangFix}) on the mathematics of finite elements:
\emph{the penalty method was the main theme of Courant's remarkable lecture; the finite element method was an afterthought.}

Today,  the penalty method is widely used in finite elements, especially in applied engineering. From a mathematical point of view, the method is non-conforming and the central question is how to choose the penalty parameter. The engineering literature is rather vague on this issue, cf. \cite{Laursen, Wriggers}. The same applies for some commercial codes that we checked, with the   COMSOL as an exception. 

 Below we will discuss this matter and we will also recall Nitsche's method which gives a simple fix to the method.

\section{The penalty method}
All the questions we are going to address can be answered by considering "mortaring" of the Poisson equation. To consider elasticity only adds more notation. To consider a contact problem where the contact zone is unknown leads to technicalities involving variational inequalities, cf. \cite{GSVel}. By treating both of these would only blur the message of this paper.

Hence, we consider the Poisson problem in a domain $\Omega\subset {\rm I\!R}^n$, that we split in two $\Omega=\Omega_1\cup \Omega_2$, with the continuity conditions at the common boundary $\Gamma=\partial \Omega_1 \cap \partial \Omega_2$:
find 
  $u=(u_1,u_2)$ such that 
\begin{equation}\label{strong} 
\begin{aligned}
    -\Delta u_i &= f \quad && \text{in $\Omega_i$,} \\
    u_1 - u_2 &= 0 && \text{on $\Gamma$,} \\
    \frac{\partial u_1}{\partial n_1} +  \frac{\partial u_2}{\partial n_2} &= 0 && \text{on $\Gamma$}, \\
    u_i &= 0 && \text{on $\partial \Omega_i \setminus \Gamma$.}
\end{aligned}
\end{equation}

\begin{figure}[h!]
   \centering
   \begin{tikzpicture}[scale=0.7]
       \draw[fill=gray!10] (0.5, 0.5) rectangle (4, 4) node[pos=.5] {$\Omega_1$};
       \draw[fill=gray!20] (4, -1) rectangle (8, 5.5) node[pos=.5] {$\Omega_2$};
       \draw (4, 0.5) -- (4, 4) node[midway, anchor=west] {$\Gamma$};
       \draw[line width=0.5mm, dashed] (4, 0.5) -- (4, 4);
      \draw[-latex] (4, 1.5) -- (5, 1.5) node[pos=0.95, anchor=north] {$\boldsymbol{n}_1$};
      \draw[-latex] (4, 3) -- (3, 3) node[pos=0.95, anchor=north] {$\boldsymbol{n}_2$};
   \end{tikzpicture}
\end{figure}

Clearly $U$, with $ U\vert_{\Omega_i}=u_i, \ i=1,2,$ solves the equation
\begin{equation}
-\Delta U=f \  \text{ in } \Omega,  \quad   \ U=0 \  \text{ on } \partial \Omega.
\end{equation}

\medskip

\emph{The continuous penalized problem} is:
Find $u_\varepsilon=(u_1^\varepsilon, u_2^\varepsilon)$ minimizing 
\begin{equation}\label{penalty}
E(v)=\sum_{i=1}^2\Big(\frac{1}{2}  \int_{\Omega_i} \vert \nabla v_i\vert ^2\, dx  -\int_{\Omega_i}  f v_i\, dx\Big)+ \int_\Gamma \frac{1} {2\varepsilon} (v_1-v_2)^2 \, ds. 
\end{equation}

It is clear that a necessary condition in the limit $\varepsilon    \to 0$ is that $u_1^0=u_2^0$ on $\Gamma$, with $\lim_{\varepsilon \to 0}u_i^{\varepsilon}=u_i^0 = u_i$, $i=1,2$. 

The physical  interpretation of non-conformity of the method is that the problem is changed to that where a distributed spring is introduced on $\Gamma$, giving 
  the interface condition
$$
 \frac{\partial u^\varepsilon_{i}}{\partial n_i} +\varepsilon^{-1} (u_{i}^\varepsilon -u ^\varepsilon_{j})=0,  \quad \quad i\not= j .
$$
\
\medskip
In the FEM we perform the minimization in subspace $V_h$  = piecewise polynomials of degree $p$.

\medskip
\emph{The Finite Element Method:} Find $ u_h^\varepsilon=(u_{1,h}^\varepsilon, u_{2,h}^\varepsilon)\in V_h=V_{1,h} \times V_{2,h}$ that minimizes the energy $ E(v)$ in the subspace $V_h$.
\medskip

The \emph{accuracy}
of the FE solution is most naturally measured in\emph{ Energy norm:}

\begin{equation}
\Vert v\Vert_E^2 = \sum_{i=1}^2 \int_{\Omega_i} \vert \nabla v_i\vert ^2\, dx + \int_\Gamma \frac{1}{\varepsilon} (v_1-v_2)^2 \, ds. 
\end{equation}

The basic mathematical result is that the \emph{FE solution is the best approximation in the energy norm}, viz.
\medskip
\begin{equation}
\Vert u_\varepsilon -u^\varepsilon_h \Vert_E \lesssim  \min_{v_h\in V_h} \Big\{ \sum_{i=1}^2 \Vert\nabla( u_i^\varepsilon - v_{i,h}  )\Vert_{0, \Omega_i}
+\varepsilon^{-1/2} \Vert u_i^\varepsilon - v_{i,h}  \Vert_{0, \Gamma}  \Big\}.
\end{equation}

For notational simplicity we now assume that $\varepsilon$ is a constant along $\Gamma$. 

Note that for FEM the exact solution is that of the penalized exact solution $u_\varepsilon$. To the total error we should thus add the error due to lack of consistency or, alternatively stated, modelling. It is not too difficult to show that this is
\begin{equation}\label{modelling}
\Vert u-u_\varepsilon\Vert_E \lesssim \varepsilon^{1/2}\Big \Vert \frac{\partial u}{\partial n} \Big \Vert_{0, \Gamma},
\end{equation}
where $ \lesssim \cdots $ stands for: there exists a positive constant $C$, such that $\leq C \cdots $.
 
 By the triangle inequality we have
\begin{equation}\label{total}
\Vert u-u^\varepsilon_h\Vert_E\leq \Vert u-u_\varepsilon\Vert_E+\Vert u_\varepsilon-u^\varepsilon _h\Vert_E.
\end{equation}
As usual, let $h$ denote the global mesh lenght. 

With the assumed smoothness of the exact penalized solution $u_\varepsilon \in H^s(\Omega)$, with $1\leq s \leq p+1$, interpolation theory yields
\begin{equation}\label{intest}
 \min_{v_h\in V_h} \Big\{ \sum_{i=1}^2 \Vert\nabla( u_i^\varepsilon - v_{i,h}  )\Vert_{0, \Omega_i}
+\varepsilon^{-1/2} \Vert u_i^\varepsilon - v_{i,h}  \Vert_{0, \Gamma}  \Big\}\lesssim  \big( h^{s-1} + \varepsilon^{-1/2}  h^{s-1/2}\big) \Vert u_\varepsilon \Vert_s.
\end{equation}
(This estimate was first given by Babu\v{s}ka \cite{Babuska}.)
By choosing 
 $\varepsilon \approx h$, we get the optimal convergence rate
\begin{equation}
\Vert u_\varepsilon-u^\varepsilon_h\Vert_E \lesssim  h^{s} \Vert u_\varepsilon \Vert_s.
\end{equation}

From \eqref{modelling} and \eqref{total} we thus get

\begin{equation}
\Vert u-u^\varepsilon _h\Vert_E \lesssim h^{1/2} \Big \Vert \frac{\partial u}{\partial n} \Big \Vert_{0, \Gamma} + h^{s} \Vert u_\varepsilon \Vert_s.
\end{equation}
From above we conclude that the modelling error is dominating, except for a very singular solution.

The modelling error can of course be made smaller by decreasing $\varepsilon $. However, \eqref{intest} shows that for a given smoothness, the interpolation error will increase.
The condition number also increases from  $\mathcal{O} (h^{-2})$, which is the one expected for a second order elliptic equation.

In conclusion: {\it with the penalty method we will never obtain a method with optimal convergence rate.}

Next, let us turn to the \emph{a posteriori error }analysis. In addition to the normal terms, and edge $E\subset \Gamma $ yields the term \cite{JuSt} 
 \begin{equation}
 \eta_E^2 =h_E+ h_E \Big\Vert \frac{\partial u^\varepsilon _{i,h}}{\partial n_i} +\varepsilon^{-1} (u^\varepsilon_{i,h} -u^\varepsilon_{j,h}) \Big\Vert_{0,E}^2,   \quad j\not=i .
 \end{equation}
 In \cite{JuSt} we also showed that the second term above is efficient, i.e. bounded by the real error if and only if $\varepsilon \approx h$, giving an additional justification for this choice.
 
 In the next section we will show that there is a very simple correction that leads to an optimal method, i.e. Nitsche's method.

\section{Nitsche's method \cite{Nitsche}}

Let us recall that the weak form obtained by minimizing \eqref{penalty} is: find $u_\varepsilon=(u_1^\varepsilon, u_2^\varepsilon) $ such that  
 \begin{equation}\label{penvar}
 \sum_{i=1}^2\Big(  \int_{\Omega_i}  \nabla u^\varepsilon _i\cdot \nabla v_i\, dx  + \int_\Gamma \frac{1 } {\varepsilon} (u^\varepsilon_i-u_j)v_{i} \, ds\Big)=
  \sum_{i=1}^2 \int_{\Omega_i}  f v_i\, dx
\end{equation}
for all test functions $v=(v_1,v_2)$.

For the solution $u=(u_1,u_2) $ of the exact problem \eqref{strong}
 it holds
\begin{equation}u_i=u_j  \quad \mbox{ on } \Gamma,
\end{equation}
and hence
 \begin{equation}
  \int_\Gamma \frac{1 } {\varepsilon} (u_i-u_j)v_{i} \, ds=0.  
\end{equation}
Further, multiplying the differential equation in \eqref{strong}, integrating over $\Omega_i,$ and integrating by parts yields
 \begin{equation}
   \int_{\Omega_i}  \nabla u_i \cdot \nabla v_i \, dx -\int_\Gamma \frac{\partial u_i}{\partial n_i} v_i  \, ds=
   \int_{\Omega_i}  f v_i\, dx . 
\end{equation}
These two relations show that the exact solution satisfies
 \begin{equation}
 \begin{aligned}
 &\sum_{i=1}^2 \Big( \int_{\Omega_i}  \nabla u_i\cdot \nabla v_i\, dx  + \int_\Gamma \frac{1} {\varepsilon} (u_i-u_j)v_{i} \, ds\Big)
 \\&-\int_\Gamma \frac{\partial u_1}{\partial n_1} v_1\, ds -\int_\Gamma \frac{\partial u_2}{\partial n_2} v_2\, ds=
 \sum_{i=1}^2  \int_{\Omega_i}  f v_i\, dx,
   \end{aligned}
\end{equation}
for all test functions $(v_1,v_2)$.

Comparing this with \eqref{penvar}, shows the non conformity of the penalty method, i.e. the exact solution does not satisfy the penalty variational form.

However, the above form  
 would already be the rudimentary Nitsche.

The next step is to use the relation
\begin{equation}
  \frac{\partial u_2}{\partial n_2} = -  \frac{\partial u_1}{\partial n_1} ,
\end{equation}
to get the weak form 
\begin{equation}
 \begin{aligned}
 &\sum_{i=1}^2  \Big(\int_{\Omega_i}  \nabla u_i\cdot \nabla v_i\, dx  + \int_\Gamma \frac{1} {\varepsilon} (u_i-u_j)v_{i} \, ds\Big)
 \\&-\int_\Gamma \frac{\partial u_1}{\partial n_1} (v_1-v_2)\, ds =
\sum_{i=1}^2   \int_{\Omega_i}  f v_i\, dx . 
   \end{aligned}
\end{equation}

However, this is not entirely satisfying, since we end up with a non symmetric method for a problem that is symmetric.

Hence, the final step is to symmetrise the problem.
For the exact solution $u=(u_1,u_2)$ it holds $u_1=u_2$ on $\Gamma$, and hence it holds
\begin{equation}
\int_\Gamma \frac{\partial v_1}{\partial n_1} (u_1-u_2)\, ds=0.
\end{equation}
We can thus add this to the variational form, and we conclude that the exact solution $u=(u_1,u_2)$ satisfies 
 \begin{equation}
 \begin{aligned}
 &\sum_{i=1}^2  \Big(\int_{\Omega_i}  \nabla u_i\cdot \nabla v_i\, dx  + \int_\Gamma \frac{1} { \varepsilon} (u_i-u_j)v_{i} \, ds\Big)
 \\&-\int_\Gamma \frac{\partial u_1}{\partial n_1} (v_1-v_2)\, ds -\int_\Gamma \frac{\partial v_1}{\partial n_1} (u_1-u_2)\, ds=
\sum_{i=1}^2   \int_{\Omega_i}  f v_i\, dx 
   \end{aligned}
\end{equation}
for all test functions $v=(v_1,v_2)$.

By denoting \begin{equation}\label{nit1}
\begin{aligned}
\B(u,v) &=
\sum_{i=1}^2  \Big(\int_{\Omega_i}  \nabla u_i\cdot \nabla v_i\, dx  + \int_\Gamma \frac{1} { \varepsilon} (u_i-u_j)v_{i} \, ds\Big)
 \\&-\int_\Gamma \frac{\partial u_1}{\partial n_1} (v_1-v_2)\, ds -\int_\Gamma \frac{\partial v_1}{\partial n_1} (u_1-u_2)\, ds,
 \end{aligned}
\end{equation}
we define the  \emph{Nitsche method}: find $u_h=(u_{1,h},u_{2,h})\in V_h$ such that 
\begin{equation}\label{nit2}
\B(u_h,v) =\sum_{i=1}^2   \int_{\Omega_i}  f v_i\, dx  \quad \forall v\in V_h.
\end{equation}
This is the one-sided Nitsche, a master-slave method, with the domain $\Omega_1$ being the slave.   

By construction, the method is \emph{consistent}. It remains to choose $\varepsilon$ so that the method is stable and have an optimal order of convergence. To this end we first have
\begin{equation}\label{coer}
\B(v,v) =\sum_{i=1}^2   \int_{\Omega_i}  \vert   \nabla v_i\vert ^2\, dx  -2 
\int_\Gamma \frac{\partial v_1}{\partial n_1} (v_1-v_2)\, ds+ \int_\Gamma \frac{1} { \varepsilon} ( v_1-v_2)^2 \, ds.
 \end{equation} 
 Schwarz inequality gives
 \begin{equation}\label{SI}
 \begin{aligned}
 &\Big\vert \int_\Gamma \frac{\partial v_1}{\partial n_1} (v_1-v_2)\, ds \Big\vert =
 \Big\vert \int_\Gamma \varepsilon^{1/2}\frac{\partial v_1}{\partial n_1} \  \varepsilon^{-1/2} (v_1-v_2)\, ds \Big\vert 
 \\
& \leq 
  \Big(\int_\Gamma\varepsilon \, \Big \vert \frac{\partial v_1}{\partial n_1}\Big  \vert^2 \, ds\Big)^{1/2}
 \Big(\int_\Gamma \frac{1} { \varepsilon} ( v_1-v_2)^2 \, ds\Big)^{1/2}.
 \end{aligned}
 \end{equation} 
 Using the inequality $  ab\leq  a^2+ \frac{1}{4} b^2, \ a,b \in \mathbb{R}$, we estimate
 \begin{equation}\label{AMG}
  \Big(\int_\Gamma\varepsilon \, \Big \vert \frac{\partial v_1}{\partial n_1}\Big  \vert^2 \, ds\Big)^{1/2}
 \Big(\int_\Gamma \frac{1} { \varepsilon} ( v_1-v_2)^2 \, ds\Big)^{1/2}
 \leq
 \int_\Gamma\varepsilon \, \Big \vert \frac{\partial v_1}{\partial n_1}\Big  \vert^2 \, ds +\frac{1}{4}
 \int_\Gamma \frac{1} { \varepsilon} ( v_1-v_2)^2 \, ds.
 \end{equation} 
 Combining \eqref{SI} and \eqref{AMG}, then gives
\begin{equation}
-2 
\int_\Gamma \frac{\partial v_1}{\partial n_1} (v_1-v_2)\, ds\geq
-2 \int_\Gamma\varepsilon \, \Big \vert \frac{\partial v_1}{\partial n_1}\Big  \vert^2 \, ds -
 \int_\Gamma \frac{1} {2 \varepsilon} ( v_1-v_2)^2 \, ds,
\end{equation}
which together with \eqref{coer} gives
\begin{equation}
\B(v,v) \geq \sum_{i=1}^2  \int_{\Omega_i}  \vert   \nabla v_i\vert ^2\, dx -2 \int_\Gamma\varepsilon \, \Big \vert \frac{\partial v_1}{\partial n_1}\Big  \vert^2 \, ds+ \int_\Gamma \frac{1} {2 \varepsilon} ( v_1-v_2)^2 \, ds.
 \end{equation}

Now, let $K\subset \Omega_1$ be an element with a side/face $E \subset \Gamma$. By a scaling argument there exists a positive constant $C_I$ such that 
\begin{equation}\label{trace}
C_I\int_K \vert   \nabla v_1\vert ^2\, dx \geq  h_E  \int_E \, \Big \vert \frac{\partial v_1}{\partial n_1}\Big  \vert^2 \, ds,
\end{equation}
where $h_E$ is the local mesh length at $E$. 
This gives 
\begin{equation}
  \int_K \vert   \nabla v_1\vert ^2\, dx - 2 \int_E \, \varepsilon \Big \vert \frac{\partial v_1}{\partial n_1}\Big  \vert^2 \, ds
   \geq \Big(  1-2\varepsilon \frac{ C_I}{h_E}   \Big) \int_K \vert   \nabla v_1\vert ^2\, dx  . 
  \end{equation}

Hence, if we choose $\varepsilon$ such that $\varepsilon\vert_E=h_E/2   \gamma$, with $\gamma >C_I$, it holds
\begin{equation}
 \int_K \vert   \nabla v_1\vert ^2\, dx - 2 \int_E \, \varepsilon \Big \vert \frac{\partial v_1}{\partial n_1}\Big  \vert^2 \, ds
  \gtrsim \int_K \vert   \nabla v_1\vert ^2\, dx. 
\end{equation}
Combining with \eqref{coer}, then yields the result.
There exists a constant $C_I$, so that the choice $\varepsilon \vert_E =h_E/2   \gamma$, with $ \gamma > C_I$, for all edges/faces on $\Gamma $, yields the {\em stability estimate}  

\begin{equation}
\B(v,v) \gtrsim \sum_{i=1}^2  \int_{\Omega_i}  \vert   \nabla v_i\vert ^2\, dx + \sum_{E\subset \Gamma} \int_\Gamma \frac{1} { h_E} ( v_1-v_2)^2 \, ds=\Vert v\Vert_E^2
 \end{equation}

 By construction, the  method is consistent, and together with the stability, one obtains an optimal \emph{a priori error estimate}:
\begin{equation}
\Vert u-u_h\Vert_E \lesssim \inf_{v\in V_h} \Vert u-v\Vert_E\lesssim h^{s-1} \Vert u\Vert_s, \quad 1\leq s\leq p+1.
\end{equation}

\medskip

\emph{Remark. } When implementing the method, the crucial question is of course how difficult it is to estimate the constant $C_I$. For piecewise linear elements this is simple. On the element $K$, with edge/face $E$,  $\nabla v_1$ is a constant vector. Hence, it holds
\begin{equation}
 \int_K \vert   \nabla v_1\vert ^2\, dx=  \vert   \nabla v_1\vert ^2  \int_K   dx  \ \mbox{ and } 
  \int_E \,  \Big \vert \frac{\partial v_1}{\partial n_1}\Big  \vert^2 \, ds =  \Big \vert \frac{\partial v_1}{\partial n_1}\Big  \vert^2   \int_E ds.
\end{equation}
Since
\begin{equation}
 \vert   \nabla v_1\vert ^2 \geq \Big \vert \frac{\partial v_1}{\partial n_1}\Big  \vert^2,
\end{equation}
it holds
\begin{equation}
 \int_K \vert   \nabla v_1\vert ^2\, dx \geq   \frac{\int_Kdx}{\int_E ds}  \int_E \, \Big \vert \frac{\partial v_1}{\partial n_1}\Big  \vert^2 \, ds.
\end{equation}
Since $\int_K dx$ is the volume/area of $K$, and $\int E ds$ is the area/length of $E$, respectively, the ratio
$ \int_Kdx/\int_E ds  $ is proportional to the diameter $h_E$, as the elements are assumed to be shape regular.

For linear elements the threshold value for the stability is then $ \int_Kdx/2\int_E ds  $ and 
 we thus have uniform stability by choosing
\begin{equation}
\varepsilon\vert _E = \alpha\,  \frac{\int_Kdx}{2\int_E ds},
\end{equation} 
with a chosen constant $\alpha <1$. 

In  \cite{MR1986022}   the following explicit estimate for the inequality \eqref{trace} is given
\begin{equation} 
\frac{n\int_Kdx}{(p+1)(p+n)\int_E ds}\int_K \vert   \nabla v_1\vert ^2\, dx \geq   \int_E \, \Big \vert \frac{\partial v_1}{\partial n_1}\Big  \vert^2 \, ds,
\end{equation}
where $p$ is the polynomial degree and $n$ is the space dimension. $\square$
\medskip

In the \emph{a posteriori estimator} the interface gives the terms 
\begin{equation}
h_E\Big \Vert  \frac{\partial u_{1,h}}{\partial n_1} +  \frac{\partial u_{2,h}}{\partial n_2}\Big \Vert_{0,E} ^2 
\end{equation}
and
\begin{equation}
h_E  \Vert u_{1,h}-u_{2,h} \Vert ^2_{0,E},
\end{equation}
for an edge/face $E\subset \Gamma$.

The estimator can be shown to be both reliable and efficient \cite{GSVmortar}.

In the formulation \eqref{nit1}--\eqref{nit2}, the role of $u_1$ and $u_2$ can be changed. The most common is, however,  to use the average 
\begin{equation}
\begin{aligned}
 \B(u,v) =
&\sum_{i=1}^2  \Big(\int_{\Omega_i}  \nabla u_i\cdot \nabla v_i\, dx  + \int_\Gamma \frac{1} { \varepsilon} (u_i-u_j)v_{i} \, ds\Big)
 \\
 &-\int_\Gamma\frac{1}{2} \Big( \frac{\partial u_1}{\partial n_1}+ \frac{\partial u_2}{\partial n_2} \Big) (v_1-v_2)\, ds 
 \\
 &-\int_\Gamma \frac{1}{2} \ \Big( \frac{\partial v_1}{\partial n_1}+ \frac{\partial v_2}{\partial n_2} \Big)(u_1-u_2)\, ds.
 \end{aligned}
 \end{equation}
 For the case when the material parameters are different in the two domains, a convex combination of fluxes can be used, cf. \cite{GSVmortar, GSVel}.
 
 \section{Numerical verifications}
 
 In this section we will numerically confirm the theoretical results. Our first test problem has the  smooth solution 
\begin{equation}
u(x,y) = x y \sin \frac{\pi x}{2} \sin \pi y
\end{equation}
in the domains
\begin{equation}
\Omega_1 =(0,1)\times (0,1) \  \mbox{ and }\  \Omega_2 =(1,2)\times (0,1).
\end{equation}
 
\begin{center} \includegraphics[width=0.7\textwidth]{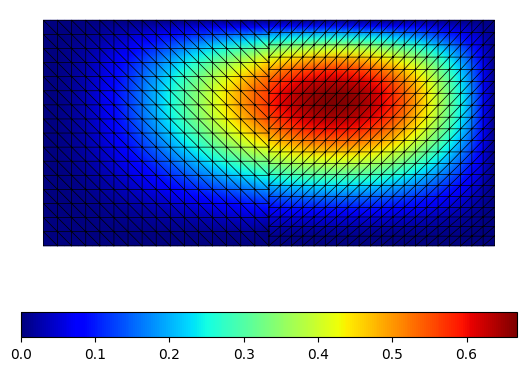}
\end{center}

\begin{center} Figure 1. The smooth test problem.
\end{center}

\bigskip

We use uniform refinements with initial meshes that do not match at the interface and linear element, $p=1$. In the penalty method we 
use $\epsilon \approx h $. Figure 2 shows the convergence with respect to the mesh length $h$. As predicted, the rates are  
$\mathcal{O}(h^{1/2})$ and  $\mathcal{O}(h)$ for the penalty and Nitsche methods, respectively.

\medskip

\begin{center}
 
\includegraphics[width=0.7\textwidth]{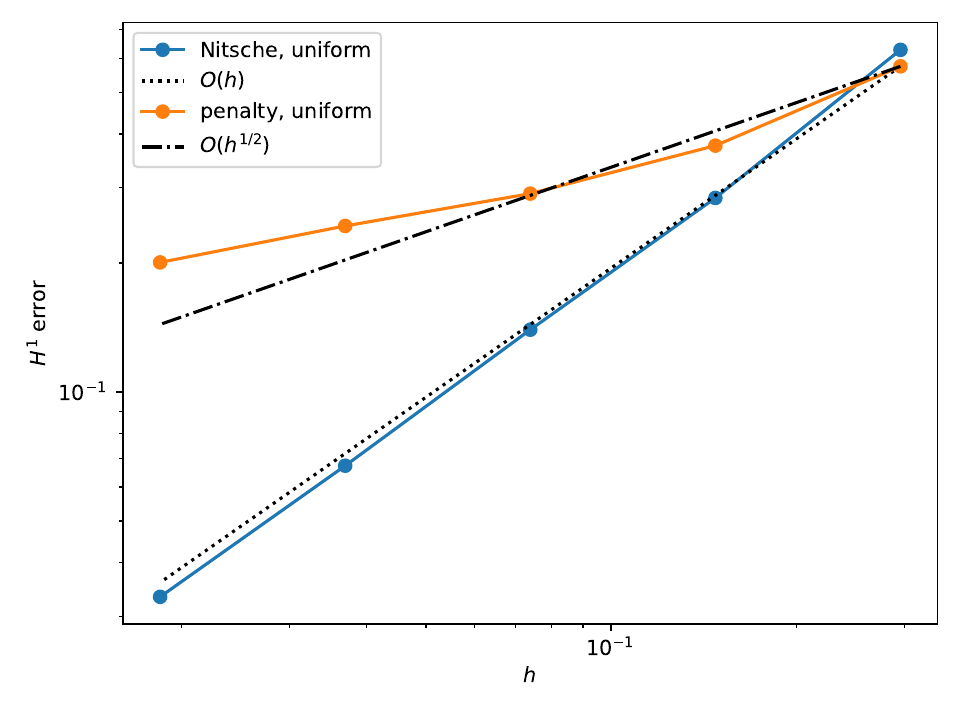}

{ Figure 2. The convergence rates for the smooth problem and a uniform mesh refinement.}
\end{center}

\bigskip
Next, we run the adaptive algorithms. We see from Figure 3, that the adaptive method detects the low convergence rate at the interface for the penalty method.

\bigskip

\begin{center}\includegraphics[width=0.7\textwidth]{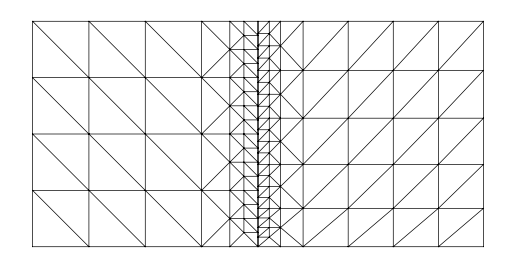}
\end{center}

\begin{center}
Figure 3. The adaptive mesh for the penalty method.
\end{center}

\bigskip

 For Nitsche any extra refinement along the interface does not occur, see Figure 4.

\begin{center}\includegraphics[width=0.7\textwidth]{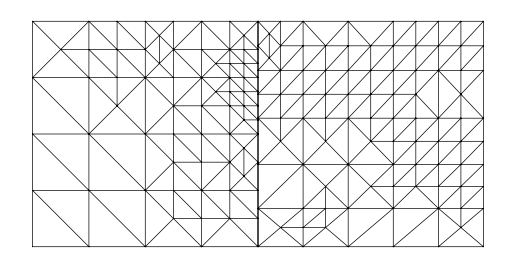}
\end{center}

\begin{center}
Figure 4. The adaptive mesh for the Nitsche method.
\end{center}

Figure 5, show that even with adaptivity the penalty falls short compared to Nitsche. Now the convergence rates are with respect to the number of degrees of freedom $N$.
\bigskip

\begin{center}

\includegraphics[width=0.8\textwidth]{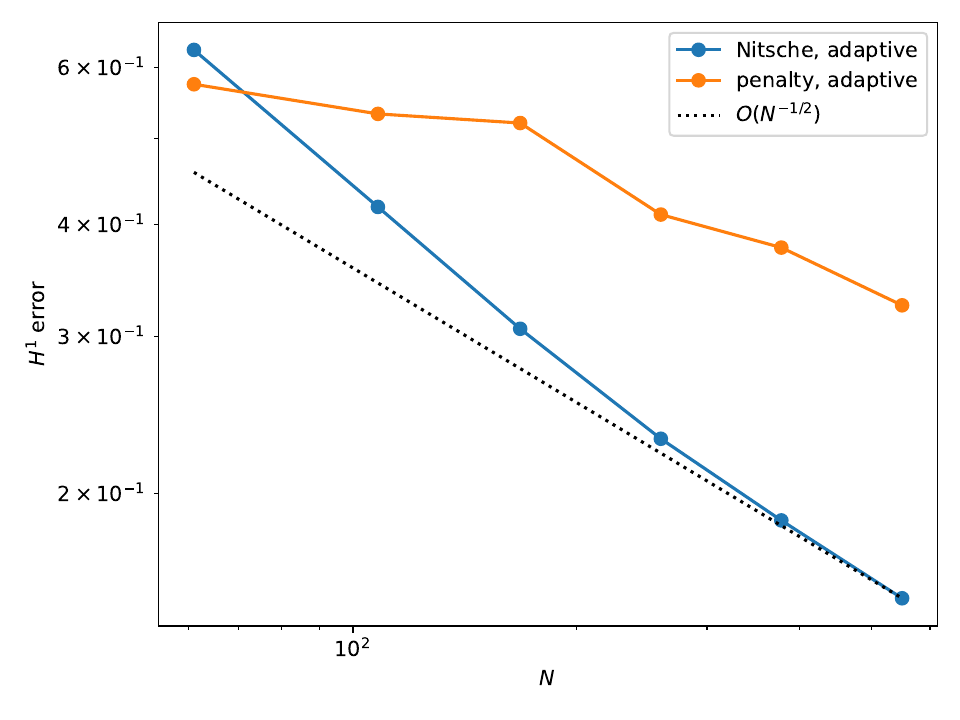}

 Figure 5. The convergence rates for the adaptive algorithms.
    \end{center}
 
 The second test example is a problem with a singularity. The domain is L-shaped and the solution is $r^{2/3}$. The solution is in $H^{5/2}$ and we choose second degree polynomials, i.e. $p=2$ in order to enforce the adaptivity. The adaptive Nitsche method detects the corner singularity and refines accordingly.

 \bigskip
  
  \begin{center}
 
\includegraphics[width=0.8\textwidth]{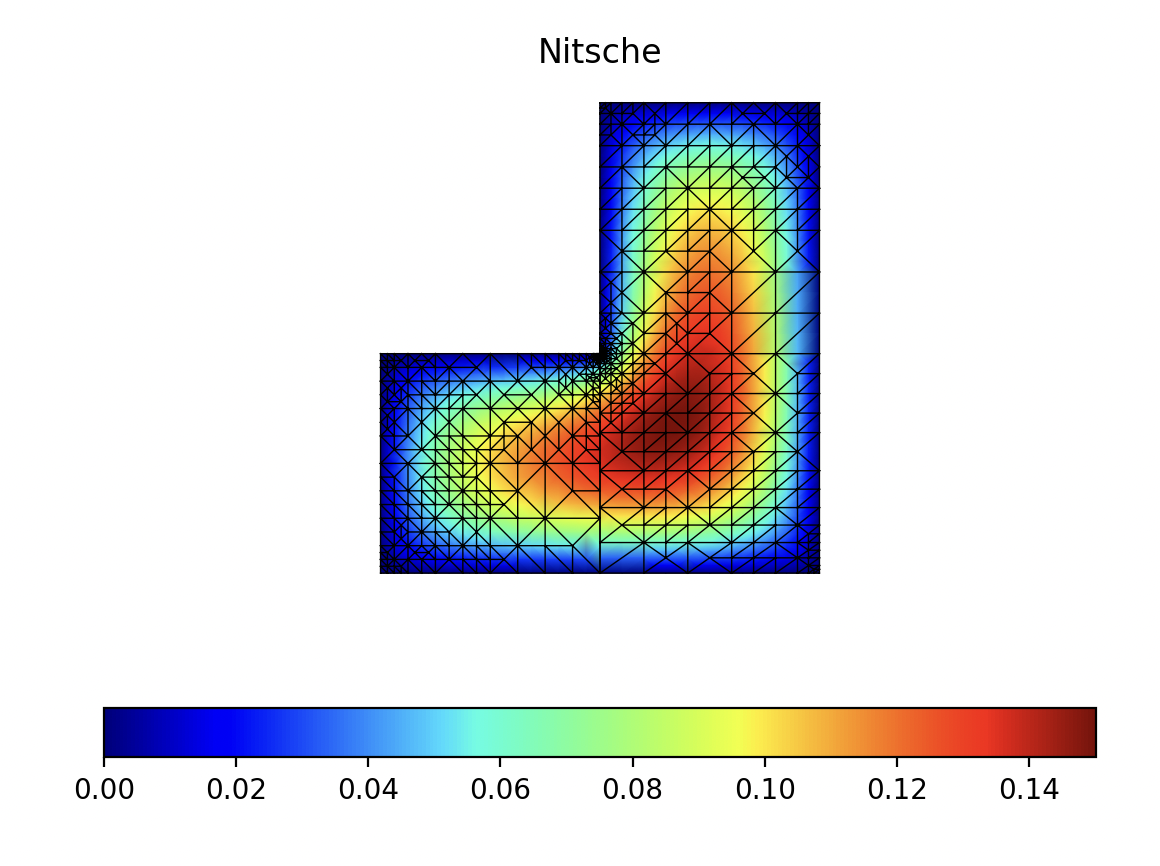} 

Figure 6. The adaptive Nitsche method for the singular test problem.

\end{center}
 \bigskip

  \bigskip
  
 \begin{center}
 
 \includegraphics[width=0.8\textwidth]{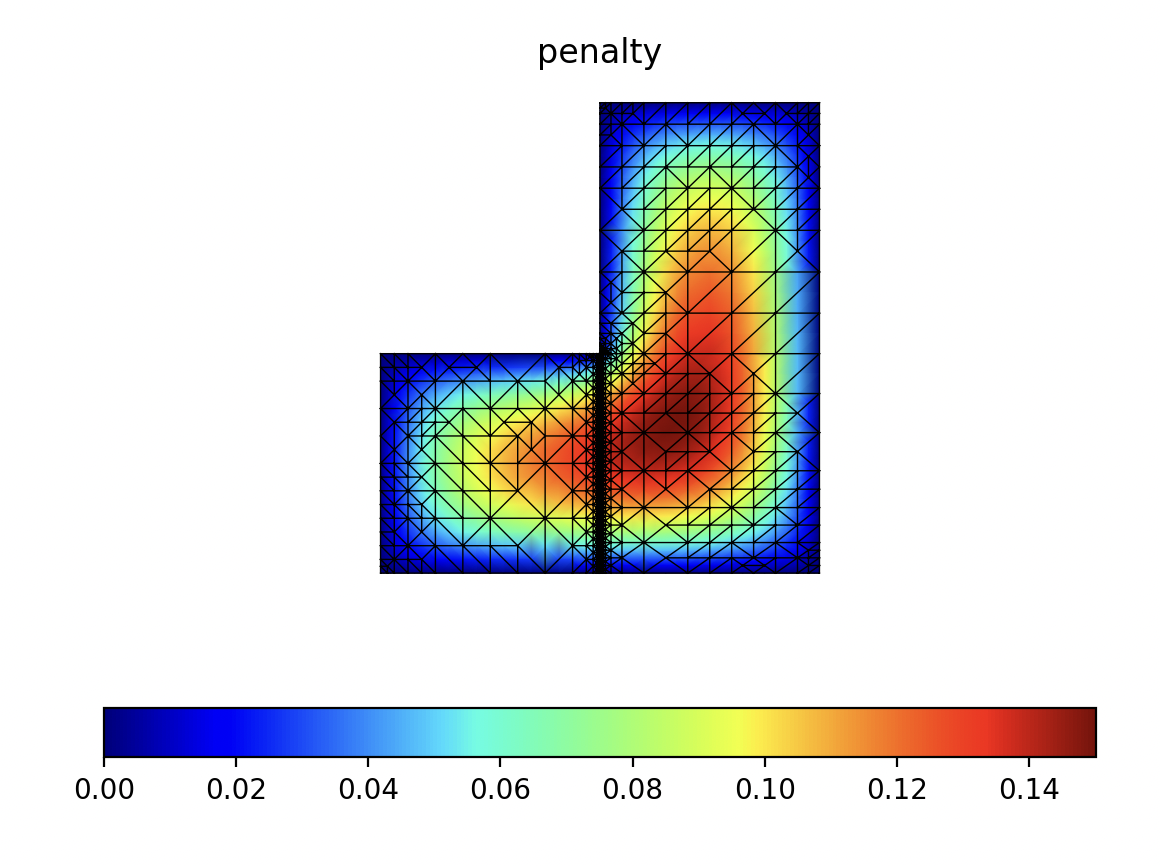}
 
 Figure 7. The adaptive penalty method for the singular test problem.
 
 \end{center}
 
 Here again, we see that the modelling error in the penalty method leads to a refinement along the whole interface.

 In the last figure the convergence rates are plotted and we again see the superiority of Nitsche's method.

 \begin{center}
 
\includegraphics[width=0.9\textwidth]{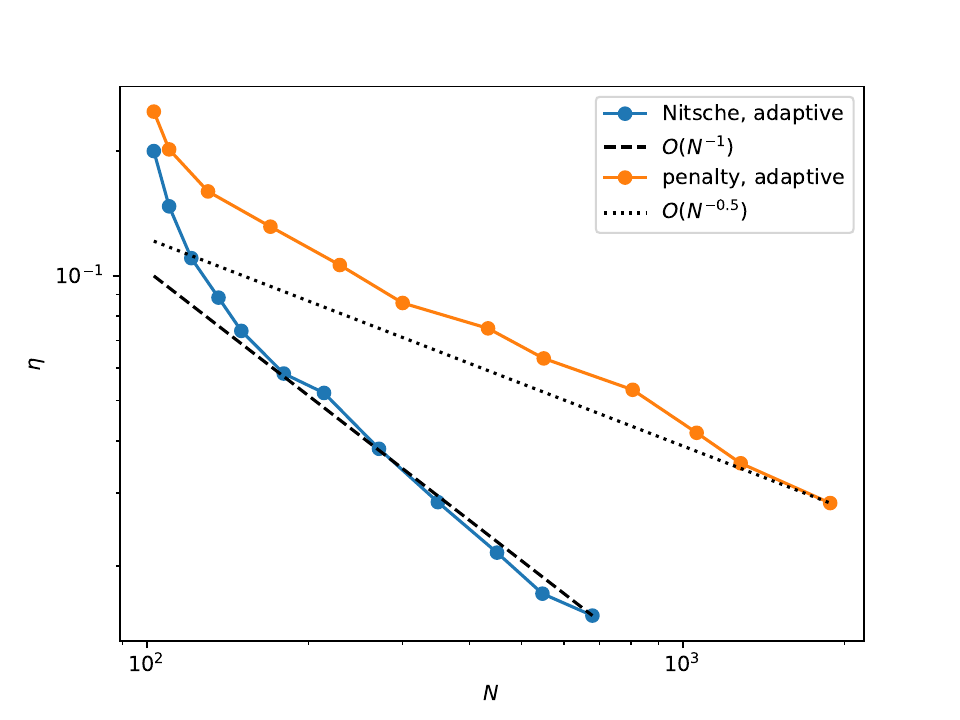}
 
 Figure 8. The convergence rates for the singular test problem. 
 \end{center}

 \section{Concluding remarks}
 
 We have shown that the Nitsche method can be seen as a very simply modification of the penalty method, yielding a method with optimal order of convergence. One can say that it is pity that Nitsche published his paper in German and in a not well-spread journal. As a consequence, it took a long time for the method to spread.
 
 However, once one understands the  method, one is stunned by its simplicity and elegance, and it is clear that there is an abundance of applications.
 
 Let us close the paper by some recollections by the senior author of this paper. He stumbled on the method by studying stabilized mixed methods for enforcing boundary conditions, and concluded that these are essentially the  same as Nitsche. At the Finnish Mechanics Days  1994, in Jyv\"askyl\"a he presented a paper \cite{RS94}. The paper was  later expanded and published internationally \cite{RS95}, and it ended with the conclusion that the method should be explored for contact problems. This was done by Jouni Freund in his thesis \cite{JF} and in a short note in the world congress of computational mechanics in Buenos Aires 1998 \cite{RSBA}. Peter Hansbo was the opponent at Freund's thesis defence, and later Hansbo and his Swedish colleagues, Mats Larson and Erik Burman, have studied the method intensively for a variety of problems. 
 
 Another group that has made substantial contributions for  elastic contact consists of Yves Renard, Patrick Hild and Franz Chouly, summarized in the recent monograph \cite{ChoulyHildRenard}. This book also include, e.g., frictional and dynamic problems.
 
 Nitsche's method is by now included in several FEA codes, like  COMSOL (\url{https://www.comsol.com})  and  GetFEM (\url{https://getfem.org}).

\bibliographystyle{plainurl}

\end{document}